\pdfoutput=1

\documentclass{amsart}
\usepackage{amsthm}
\usepackage{amsmath}
\usepackage{amssymb}
\usepackage{epsfig}
\usepackage{amsmath}
\usepackage{amsfonts}
\usepackage{amscd}
\usepackage{amssymb}
\usepackage[all]{xy}

\theoremstyle{plain}
\newtheorem{theorem}{Theorem}

\newtheorem{corollary}{Corollary}

\theoremstyle{definition}

\newtheorem{remark}[theorem]{Remark}

\newcommand\tr{\operatorname{tr}}

\newcommand\dist{\operatorname{distance}}
\newcommand\otr{\operatorname{otr}}
\newcommand\Isom{\operatorname{Isom}}
\newcommand\Relength{\operatorname{Relength}}
\newcommand\length{\operatorname{length}}
\newcommand\Arccosh{\operatorname{Arccosh}}
\newcommand\Log{\operatorname{Log}}

\begin{document}
\begin{center}
\textbf{Exceptional hyperbolic 3-manifolds}

\medskip
David Gabai, Maria Trnkova
\end{center}
\medskip

We correct and complete a conjecture of D. Gabai, R. Meyerhoff and N. Thurston on the classification and properties of thin tubed closed hyperbolic 3-manifolds.  We additionally show that if $N$ is a closed hyperbolic 3-manifold, then either $N=Vol3$ or $N$ contains a closed geodesic that is the core of an embedded tube of radius $\log(3)/2$.

\footnote{


The first author was partially supported by grants NSF DMS-0854969 and NSF DMS-1006553.

The second author was partially supported by the grant  NSF DMS-0854969 and by grant P201/11/0356 of The Czech Science Foundation.}

\section{Introduction}

An \emph{exceptional} hyperbolic 3-manifold is a closed hyperbolic 3-manifold which does not have an embedded hyperbolic tube of radius $\log (3)/2$ about each of its shortest geodesics.    They were introduced in [G] where geometric and topological rigidity theorems were proven for nonexceptional manifolds.  A detailed investigation of exceptional manifolds was conducted in [GMT] where the corresponding rigidity theorems were extended to all closed hyperbolic 3-manifolds.  Those results in turn were used in [G1] to prove the Smale conjecture for closed hyperbolic 3-manifolds, i.e. the inclusion Isom($N$)$\to$ Diff($N$) is a homotopy equivalence.   Properties of exceptional manifolds were also used in [GMM] (in conjunction with [GM]) to establish a lower bound on the volume of a closed hyperbolic 3-manifold, giving a 100+ improvement on the previously known lower bound.   They were used in the work of Agol [A] and Agol-Dunfield [AST] to improve the lower bound and give other estimates that were essentially used in [GMM] to show that the Weeks manifold is the unique closed hyperbolic 3-manifold of minimal volume.  Properties of exceptional manifolds were also used in [ACS] to give volume bounds for other classes of 3-manifolds.

An exceptional manifold $N$ gives rise to a \emph{marked} 2-generator subgroup $G$ of $\pi_1(N)$ generated by elements $f$ and $w$ where the axis $\delta_0\subset \mathbb H^3$ of $f$ projects to a shortest geodesic of $N$ and the element $w$ sends $\delta_0$ to a nearest covering translate $\delta_1$ with $d(\delta_0,\delta_1)\le \log(3).$   In [GMT] the set of marked 2-generator groups  arising from exceptional manifolds is identified with a subset $S=\exp(T)$ of a compact region  of $\mathbb C^3$.  Furthermore this region can be chopped up into about a billion regions and that any marked 2-generator group arising as above lies in one of seven small \emph{exceptional regions} $X_i$, $i=0,...,6$.  Each such region $X$ has a \textit{quasi-relator} $r(X)$, i.e.  a word in $f$, $w$, $F=f^{-1}$, $W=w^{-1}$ that is very close to the identity at all points inside the region $X$. For more detail see Chapter 1 of [GMT].  The authors of [GMT] made the following conjectures about the exceptional manifolds and the exceptional regions:

\medskip

\textbf{Conjecture.} (Exceptional Manifolds Conjecture) \label{gmt conjecture} \textit{Each exceptional box $X_i$, $0\leq i\leq6$, contains a unique element $s_i$ of $S$. Further, if $\{G_i,f_i,g_i\}$ is the marked group associated to $s_i$ then $N_i=\mathbb H^3/G_i$ is a closed hyperbolic 3-manifold with the following properties:
\begin{description}
\item[i] $N_i$ has fundamental group $<f,w;r_1(X_i), r_2(X_i)>$, where $r_1(X_i)$, $r_2(X_i)$ are the quasi-relators associated to the box $X_i$.
\item[ii] $N_i$ has a Heegaard genus-2 splitting realizing the above group presentation.
\item[iii]  $N_i$ nontrivially covers no manifold.
\item[iv] $N_6$ is isometric to $N_5$.
\item[v] If $(L_i,D_i,R_i)$ is the parameter in $T$ corresponding to $s_i$, then $L_i,D_i,R_i$ are related as follows:\\
For $X_0, X_5, X_6$: $\quad L=D, R=0$.\\
For $X_1, X_2, X_3, X_4$: $\quad R=L/2$.
\end{description}}

\medskip

It was shown in [GMT] that for each $j$ there is a Heegaard genus-2 manifold $M_j$ with presentation as in (i).  Also the closed hyperbolic 3-manifold Vol3 is the unique exceptional manifold $N_0$ corresponding to the region $X_0$ and (v) holds for $X_0$.

K. Jones and A. Reid [JR]  proved that $N_0$ nontrivially covers no manifold and found arithmetic hyperbolic manifolds $N_i$ for all the regions $X_i$, $i=0,1,2,4,5,6$.  They also showed that the exceptional manifolds $N_5$ and $N_6$ are isometric.

The manifold $N_3$ is described by M. Lipyanskiy experimentally in [L].

A. Champanerkar, J. Lewis, M Lipyanskiy and S. Meltzer [CLLMR] proved that each exceptional region contains a unique hyperbolic 3-manifold $N_i$.    They also established  properties $(i)$ and $(ii)$ for all the $N_i$'s and $(v)$ for all the $X_i$'s.   As each $N_i$ is a rational homology 3-sphere, they observe that no $N_i$ can cover a non orientable 3-manifold.

A. Reid [CLLMR] proved that the manifolds $N_1$ and $N_5=N_6$ nontrivially cover no manifold.  \vskip10pt

\noindent\textbf{What needs to be done.} To complete the proof of the Exceptional Manifolds conjecture, it suffices to show that the manifolds  $N_2$, $N_3$ and $N_4$ nontrivially cover no orientable manifold. The main result of this paper is a positive proof of the conjecture iii) reduced by the following results.
\medskip

\begin{theorem} \label{finish}  If $N_i$ is an exceptional manifold and $p:N_i\to M$ is a nontrivial covering projection, then up to conjugacy either $p:N_2\to m010(-2,3)$ or $p:N_4\to m371(1,3)$.  Furthermore, any two such coverings (for a given domain and range) are topologically conjugate.  Finally  for each $p$, $\deg(p)=2.$   \end{theorem} \vskip10pt


\begin{remark}  i)  For a given $N_i, i=2,4$ there are three different homotopy classes of covering projections as above.

ii)  The manifolds $m010(-2,3)$ and $m371(1,3)$ are given by SnapPea notation. \end{remark}


\medskip

\begin{remark} \label{quotients} Neither of $m010(-2,3)$, $m371(1,3)$ are exceptional manifolds as their shortest geodesics have tube radii $> \log(3)/2$.\end{remark}

\begin{corollary}  Any exceptional manifold is isometric to one of $N_0, N_1, N_2, N_3, N_4$ or $N_5$.
\end{corollary}

In the course of proving Theorem \ref{finish} we obtain the following.

\begin{theorem}\label{thick tubes} If $N_i\neq N_0$, then some geodesic in $N_i$ is the core of an embedded tube of radius log(3)/2.  \end{theorem}

\begin{corollary} \label{vol3} $Vol3$ is the unique closed hyperbolic 3-manifold such that no closed geodesic is the core of an embedded tube of radius log(3)/2.
\end{corollary}

\medskip

The proofs of Theorem \ref{finish} and \ref{thick tubes} require rigorous computer assistance.  They are motivated by output from the computer programs Snap [GHN] and SnapPea [W].  The program Snap [GHN] studies arithmetic and numerical invariants of hyperbolic 3-manifolds and is based on the program SnapPea [W] and on the number theory package Pari. SnapPea was written by Jeff Weeks' for studying hyperbolic 3-manifolds and Pari calculates arithmetic and number theoretic functions with high precision.

In Section 2 we use the length and ortholength information provided by Snap for the manifolds $N_i$, $i=2,3,4$ to conclude that $N_3$ nontrivally covers no manifold and that  $N_2$, $N_4$ can only nontrivally cover a manifold via a special 2-fold one.  Rigorously verifying output of SnapPea we then show that these manifolds actually have 2-fold quotients and they are exactly as in Theorem \ref{finish}.

In order to make our results rigorous we have to check data from Snap.  Snap makes all calculations with high precision and works well in practice but some algorithms (like algorithms for length and ortholength spectra) depend on the limitations of fixed-precision floating point computations.   Rigorous evaluations of the length and ortholength spectra for the manifolds $N_i$ are also needed to prove Theorems \ref{thick tubes}  and Remark \ref{quotients}.   The second author wrote a package in Mathematica to rigorously compute length and ortholength spectra for manifolds. This program can be used independently of this paper. The theoretical part of the algorithm, together with an explanation of its rigor, is given in section 3.
\vskip 8pt
\noindent\textbf{Acknowledgments:}  M. Trnkova is very thankful to J. Weeks, N. Dunfield and C. Hodgson for their help with computer program Snap, the Mathematical Department of Princeton University for hospitality and the Institute for Advanced Study for using their computer cluster. She was a visiting student research collaborator of D. Gabai at Princeton University during the preparation of this manuscript and reported on these results at the NSF sponsored March 14-16, 2011 FRG conference in  Princeton.


\section{Ortholines of thick tube geodesics}

This chapter completes the proof of the theorem assuming the correctness of various results of Snap up to 50 decimal places. In the next chapter  we rigorously check the needed data.

Recall the following three facts from [CLLMR].  For each exceptional family there is a unique manifold $N_i$, $i=2,3,4$. Closed hyperbolic 3-manifolds $s778(-3,1)$ and $v2018(2,1)$ from SnapPea's census are isometric to $N_2$. There are no manifolds in SnapPea's census isometric to $N_3$ or $N_4$ because of their large volumes.

Matrix representatives for the generators $f, w$ of fundamental group $G_i$ depends only on the three complex parameters $L',D',R'$ from the exceptional region $X_i$ [GMT]:

\begin{equation}\label{generator1}
f=
\left(
  \begin{array}{cc}
    \sqrt{L'} & 0 \\
    0 & {1}/{\sqrt{L'}} \\
  \end{array}
\right)
\end{equation}

\begin{equation}\label{generator2}
w=\left(
  \begin{array}{cc}
    \frac{\sqrt{R'}\dot (\sqrt{D'}+1/\sqrt{D'})}{2} & \frac{\sqrt{R'}\dot (\sqrt{D'}-1/\sqrt{D'})}{2} \\
    \frac{(\sqrt{D'}-1/\sqrt{D'})}{2\sqrt{R'}} & \frac{(\sqrt{D'}+1/\sqrt{D'})}{2\sqrt{R'}} \\
  \end{array}
\right)
\end{equation}
For example, the region $X_2$ is approximately bounded by:
$$
\begin{array}{llll}
 l'_{min}= -1.78701& l'_{max}=-1.78527  & t'_{min}=-2.27253 & t'_{max}=-2.27130  \\
 d'_{min}=-1.07428 & d'_{max}=-1.07273  & b'_{min}=-2.71846 & b'_{max}=-2.71736  \\
 r'_{min}=0.74163 & r'_{max}=0.74301  &  a'_{min}=-1.52929 & a'_{max}=-1.52832
\end{array}
$$
The quasi-relators $r_1(X_i),r_2(X_i)$ are relators of $G_i$ at some triple $(L'_i,D'_i,R'_i)\in X_i$, $L'=l'+t'$, $D'=d'+b'$, $R'=r'+a'$.
K. Jones and A. Reid [JR] showed that a two generator group can be determined up to conjugacy by the triple of traces $(\tr f^2,\tr w^2,\tr f^2w^2)$. This triple generates a number field which is the invariant trace field.


[CLLMR] computed the trace triple $p=\tr f$, $q=\tr w$, $r=\tr f^{-1}w$ up to high precision for all exceptional manifolds, e.g. 100 significant digits, and represented them as roots of polynomials over integers. Therefore we can get generators $f, w$ with high precision as well. The initial triple and trace triple are connected by equations:

\begin{equation*}
L'=\left( \frac{p\pm \sqrt{p^2-4}}{4} \right)^2
\end{equation*}
\begin{equation*}
D'=\left(\frac{q\sqrt{R'}\pm \sqrt{q^2R'-(1+R')^2}}{1+R'}\right)^2
\end{equation*}
\begin{equation*}
R'=\frac{qL'-r\sqrt{L'}}{r\sqrt{L'}-q}
\end{equation*}
After calculating entries of generators (\ref{generator1}), (\ref{generator2}) we can determine a manifold in Snap which is isometric to $N_i$. We recall to this manifolds as $N_i$. For the manifold $N_i$ Snap calculates length and ortholength spectra. Our proof is based on analyzing this information.

\subsection{Manifold N2}
For the region $X_2$ there exists a unique manifold $N_2$ isometric to $s778(-3,1), v2018(2,1)$ [CLLMR]. The fundamental group of this manifold
is:
$$
<f,w\,|\;FwfwfWffWfwfwFww,\; FFwFFwwFwfwfwFww>.
$$

The triple $(L',\,D',\,R')$ from the exceptional box $X_2\subset C^3$ is presented by roots of polynomials over integers. Therefore we can get entries of generators with any desirable precision. We display here only 5 decimal digits.
$$f=
\left(
  \begin{array}{cc}
0.74293-1.52908 i & 0+0 i \\
    0+0 i & 0.25706+0.52908 i \\
  \end{array}
\right)
$$
$$
w=\left(
  \begin{array}{cc}
    0.39135-0.96022 i & -0.30677-1.26724 i \\
    0.59162-0.48807 i & 0.60864-0.03977 i \\
  \end{array}
\right)
$$

According to [JR] and Snap the volume of $N_2$ is 3.6638... The orientable closed hyperbolic 3-manifold $m003(-3,1)$ has the smallest volume 0.9427... [GMM]. So, the manifold
$N_2$ is not a $p$-fold cover of an orientable closed hyperbolic manifold $M$ if $p>3$.

If $M$ is not an exceptional manifold, then any shortest closed geodesic of the
underlying manifold $M$ must have a tube of radius more than
$\log{3}/2$ [GMT]. Consequently, any component $\delta$ of its preimage  must be a closed geodesic
with a tuberadius $> \log{3}/2$.  Since the density of the volume of a tube $W$ about a geodesic in $N_2$ is less or equal to the 0.9 by [P] it follows that $volume(W)<3.3341$. Hence any geodesic in $N_2$ with tuberadius $> \log{3}/2$ must have length $< 3.18385$.  Recall that for a tube $W$ we have the equality $volume(W)=\pi l
\sinh^2{r}$ where $l$ is the length of the core geodesic and $r$ is the radius of the tube.




From Snap we get the list of all geodesics of length up to 3.18385 of manifold $N_2$.

\begin{table}[htdp]
\begin{center}
Length spectrum for $N_2$
\begin{tabular}{|c|c|c|c|}
\hline
 orbit & geodesic length & shortest ortholine & geodesic number\\
\hline
  0 & 1.06128-2.23704*i & 1.07253-1.94716*i & 0,1,3,6,7,8 \\
\hline
 1 & 1.06128+2.23704*i & \bf{1.52857-1.14372*i} & 2,4,5 \\
\hline
 2 &  1.76275+3.14159*i & 1.06128+2.23704*i & 9,10,11 \\
\hline
 3 &  2.13862-0.79928*i & 0.95994+1.31100*i & 12,13,14,15,16,17\\
\hline
 ... & & & \\
\hline
\end{tabular}
\end{center}
\end{table}%

Here we adapt Snap terminology. The \emph{orbits} are given by the action of Isom$(N_2)$ on $N_2$.  An \emph{ortholine} (also called \emph{orthocurve} in the literature e.g. [GMT]) is a geodesic segment which runs perpendicular from one closed geodesic to itself or from one geodesic to another.  It's \emph{length} (also called the \emph{orthodistance} in [GMT]) is a complex number whose imaginary part is well defined mod $2\pi$ if either both geodesics are oriented or the curve goes from a geodesic to itself.  Note that the \emph{tube radius} of a geodesic is one half of the real length of a shortest ortholine to itself.
Geodesics which do not appear in the table have shortest ortholines to itself of length less than $\log(3)=1.0986...$, hence have tube radius $<\log(3)/2$.  Therefore only geodesics in orbit [1] have tube radius $\ge \log(3)/2$ and these have tube radii approximately 0.764285. According to Snap, the real length of the shortest ortholine between any two distinct elements of orbit [1] is 0.88137.
This implies that if $p:N_2\to M$ is a non trivial covering projection, then it cannot happen that $p(\delta)=p(\delta')$ where $\delta,\delta'$ are distinct elements of orbit [1] and $p(\delta)$ is a shortest geodesic of $M$.


 Let $\delta$ denote an element of orbit [1].  After fixing a normal vector $z$ to $\delta$ and an orientation on $\delta$, then the initial and final points of an oriented ortholine respectively naturally determine complex numbers (also known as the \emph{basings}) well defined up to length($\delta)$.  Indeed, if the initial (resp. final) point corresponds to the number $v$ (resp. $w$), then translating $z$ distance $v$ (resp. $w$) along $\delta$ takes $z$ to the initial (resp. minus the final) tangent vector to the ortholine.  The following table lists the ortholine spectrum to $\delta$ for ortholines whose real length is at most 2.


{\footnotesize
\begin{table}[htdp]
\begin{center}
Orthospectrum of a geodesic from orbit [1] for $N_2$
\begin{tabular}{|c|c|c|}
\hline
ortholine length & initial point of ortholine & final point of ortholine\\
\hline
1.52857-1.14372*i  & -0.06128-0.99023*i &  0.46936+0.12829*i \\
\hline
1.52857-1.14372*i & -0.06128+2.15137*i &  0.46936-3.01330*i \\
\hline
1.76275+3.14159*i &  0.20404+1.13983*i & 0.20404-2.00176*i \\
\hline
1.76275+3.14159*i & -0.32660+0.02131*i & -0.32660-3.12028*i \\
\hline
1.96864+2.53545*i  & -0.06128+0.58057*i &  0.46936+1.69909*i \\
\hline
1.96864+2.53545*i & -0.06128-2.56102*i & 0.46936-1.44250*i \\
\hline
 ... & & \\
\hline
\end{tabular}
\end{center}
\end{table}%
}

The information of this table is schematically drawn in Figure 1.   Observe that the union of the first and the third pairs of ortholines form closed geodesics as the complex parts of their endpoints differ by $\pi$ while the other ortholines are closed geodesics themselves.

\begin{figure}
\includegraphics[scale=0.45]{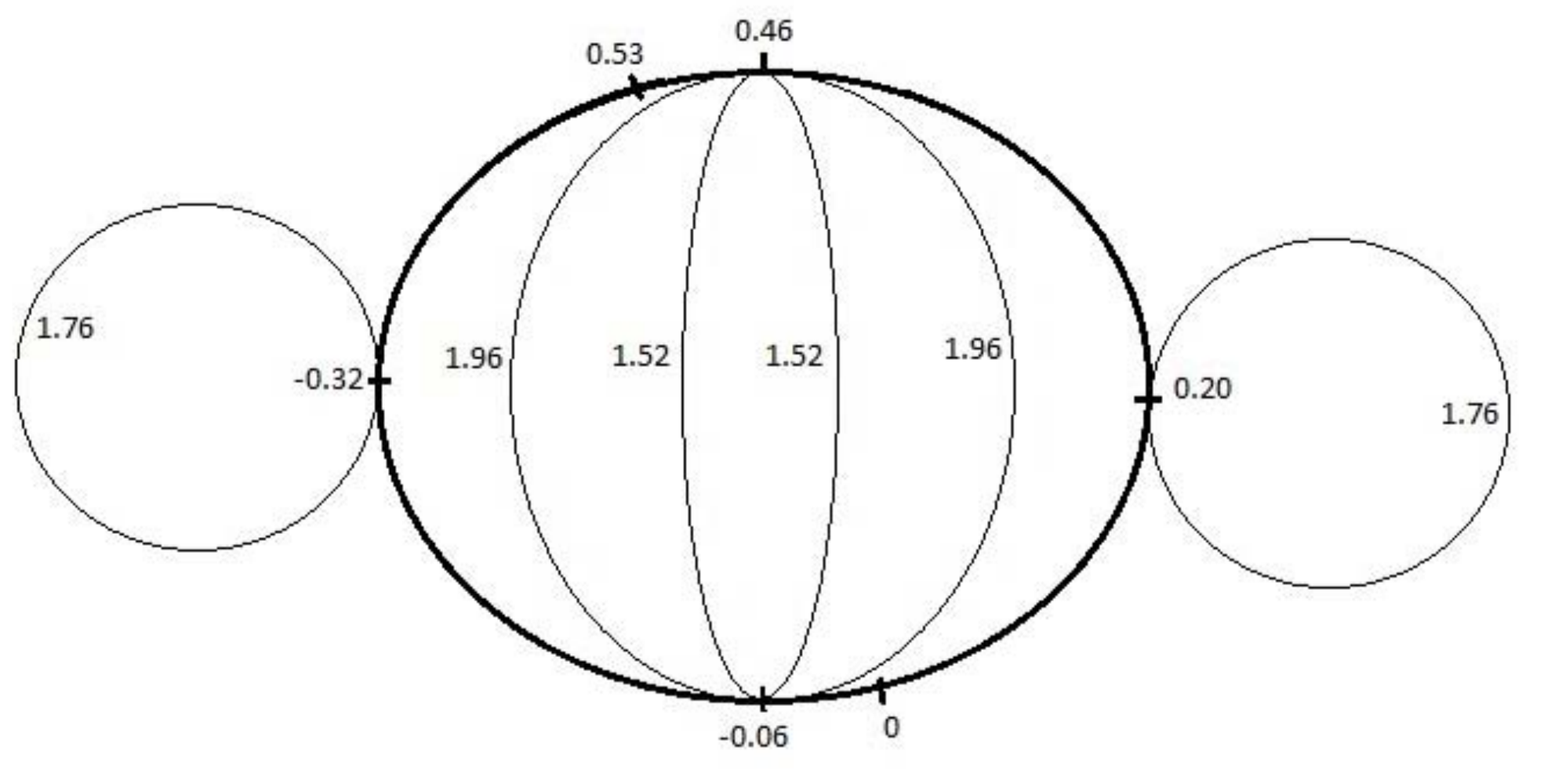}
\caption{Ortholines of geodesic from orbit [1].}
\end{figure}

It now follows that $N_2$ does not 3-fold cover any manifold $M$.  If so, then as discussed above some geodesic in orbit [1], say $\delta$, 3-fold covers a shortest geodesic $\kappa$ of $M$.  Let $\omega$ be a  shortest oriented ortholine $\omega$ from $\kappa$ to itself.  As $\omega$ is the image of an interval, it's preimage will consist of three ortholines connecting $\delta$ to itself, whose initial points  will be spaced at real distance (Re(length($\delta$))/3 along $\delta$.     Since each lift will correspond to a shortest ortholine of $\delta$ to itself, this contradicts the ortholine spectrum depicted in Figure 1.

Notice that the ortholine spectrum is not inconsistent with the existence of a 2-fold covering space.  We searched the SnapPea census for manifolds $M$ with volume$(M)$=volume$(N_2)/2$, then applied the covering space command to create 2-fold covers and then looked for manifolds whose first homology agreed with that of $N_2$.  We then discovered that $m010(-2,3)$ has a double cover  $m010\sim 0(1,0)$ with fundamental group isomorphic to that of $\pi_1(N_2)$.  Next we applied SnapPy [CD] to show that
the double cover $m010\sim 0(1,0)$ is isometric to manifolds $s778(-3,1), v2018(2,1)$. This proves the existence of a degree-2 covering map $p:N_2\to m010(-2,3)$.





Now consider any 2-fold covering $q:N_2\to M$.  The fundamental group $\pi_1(N_2)$ is a subgroup of index 2 of $\pi_1(M)$ and hence is normal which implies that all 2-fold covers arise as free $\mathbb Z_2$ symmetries of $N_2$. Looking at the orthospectrum of geodesics of orbit [1] we see that there is at most one way to map each geodesic onto itself such that the orientation is preserved. At the same time the other two geodesics are each mapped to the other. Thus there are no more than three free $\mathbb Z_2$ symmetries of $N_2$.

From Snap we know the symmetry group of $N_2$.  Inspection of the isometry group shows that three symmetries are free $\mathbb Z_2$:
\begin{enumerate}
\item $f\rightarrow WfWWffW,\quad w\rightarrow WfWFWF,$
\item $f\rightarrow fwfwFw,\quad w\rightarrow fwfWffW,$
\item $f\rightarrow wFFwFWF,\quad w\rightarrow wFFwwFw.$
\end{enumerate}
Snap uses labels 2, 4, 5 for the geodesics of orbit [1]. The first isometry (resp. second, third) preserves the geodesic number 4 (resp. 5, 2). These geodesics have the following representations in terms of generators of the fundamental group:
$$2\rightarrow WF, \quad 4\rightarrow WfW, \quad 5\rightarrow ffW.$$
Using these three isometries (1), (2), (3) we can map a geodesic of orbit [1] to any other geodesic of the same orbit. It follows that $q$ is conjugate to $p$ where the conjugating map is one of these isometries.




{\footnotesize
\begin{table}[htdp]
\begin{center}
Ortholine spectrum up to 1.4 between geodesics from orbit [1] for $N_2$. \\
For each ortholine integer numbers $2$, $4$ or $5$ denote name of geodesics \\
where the ortholine has its endpoints.

\begin{tabular}{|l|l|l|}
\hline
ortholine length & initial point of ortholine & final point of ortholine\\
\hline
0.88137-1.57080*i  & 2: 0.20404+1.13983*i & 4: 0.32660-2.40328*i \\
\hline
 0.88137-1.57080*i  & 2: 0.20404-2.00176*i & 4: 0.32660+0.73831*i \\
\hline
 0.88137-1.57080*i  & 4:-0.20404-0.38021*i & 5:-0.32660-1.18408*i \\
\hline
 0.88137-1.57080*i  & 4:-0.20404+2.76138*i & 5:-0.32660+1.95751*i \\
\hline
 0.88137+1.57080*i  & 2:-0.32660+0.02131*i & 5: 0.20404+3.07603*i \\
\hline
 0.88137+1.57080*i  & 2:-0.32660-3.12028*i & 5: 0.20404-0.06557*i \\
\hline
 ... & & \\
\hline
\end{tabular}
\end{center}
\end{table}%
}

 We move to the case $p:N_2\rightarrow M$ when $M$ is an exceptional manifold.


   The paper [GMT] proves that if $M$ is an exceptional manifold, then there is another exceptional manifold $N$ such that $N$ arises from a parameter in one of the seven exceptional boxes.  The fundamental group of $N$ is generated by two elements $f,w$ where $f$ is a hyperbolic element whose axis $B$ projects to a shortest geodesic and $w$ is a hyperbolic element which takes $B$ to a nearest translate.  In particular $N$ is a covering of $M$.

We now show that $N=N_2$. By [GMM] we know that $vol(N)>0.9427$. This would imply that the ratio $vol(N_2)/vol(N)=a/b$ where $a,b$
are integers at most 8 which is not true as we know from [CLLMR] all the 7 possibilities for $N$. Since $N=N_2$ any element $\gamma$ of orbit [0] maps with degree-1 to a shortest
geodesic $\kappa\subset M$.  Furthermore, any shortest
ortholine for $\kappa$ has the same real length of a shortest self ortholine for $\gamma$. It follows that
if $W$ is a maximal tube about $\kappa$, then $vol(W)>1.05433$. By [P] the volume of the tube
$W$ has density at most 0.91 in $M$ and hence $vol(M)>1.15861$.  Therefore $\deg(p) \leq 3$.  Thus the
preimage of $\kappa$ consists of $\gamma$ together with other geodesics which together map with
degree at most 3 to $\kappa$.

There are no geodesics in $N_2$ with real length twice that of $\gamma$, thus the preimage of $\kappa$ lies entirely in orbit [0]. The
real length of shortest ortholines between sets of geodesics \{0, 1, 3\} and \{6, 7, 8\} is 0.21561. It is less then the shortest ortholine of $\omega$ therefore there is no 3-fold cover.

As analyzed earlier the shortest geodesic in a 2-fold quotient $M$ has length approximately 0.53064 which implies that it is a non exceptional manifold.   Therefore, $N_2$ nontrivially covers no exceptional manifold if the length and ortholength spectra are rigorous.


{\footnotesize
\begin{table}[htdp]
\begin{center}
Ortholine spectrum up to 1.52 between geodesics from orbit [0] for $N_2$
\begin{tabular}{|c|c|c|}
\hline
ortholine length & initial point of ortholine & final point of ortholine\\
\hline
0.21561-1.16921*i &  0: 0.42283+1.39391*i & 1: 0.10781-1.66594*i \\
 0.21561-1.16921*i & 1:-0.42283+2.59417*i & 3:-0.10781+1.08850*i \\
 0.21561+1.97238*i & 0:-0.10781-0.62917*i & 3: 0.42283+3.11158*i \\
 0.21561+1.97238*i & 6:-0.05913-2.35087*i & 7: 0.47151-1.69595*i \\
 0.21561+1.97238*i & 6: 0.47151-0.32780*i & 8:-0.05913+0.34257*i \\
 0.21561+1.97238*i & 7:-0.05913+2.56416*i & 8: 0.47151+2.36564*i \\
 1.07253-1.94716*i & 0: 0.15751+0.38237*i & 0:-0.37312-1.64070*i \\
 1.07253-1.94716*i & 1:-0.15751-2.67748*i & 1: 0.37312-0.65440*i  \\
 1.07253-1.94716*i & 3: 0.15751+2.10004*i & 3:-0.37312+0.07697*i \\
 1.07253-1.94716*i & 6: 0.20619-1.33934*i & 6:-0.32445+2.92077*i \\
 1.07253-1.94716*i & 7: 0.20619-2.70749*i & 7:-0.32445+1.55262*i \\
 1.07253-1.94716*i & 8: 0.20619+1.35410*i & 8:-0.32445-0.66897*i \\
 ... & & \\
\hline
\end{tabular}
\end{center}
\end{table}%
}

\subsection{Manifold N3}
In a similar way we show that manifold $N_3$ cannot nontrivially cover any orientable manifold. The fundamental group of this manifold and the value of generators (up to 5 decimals) in the box $X_3$ is:
$$
<f,w\,|\;FFwfwFFwwFWFwFWfWFWffWFWfWFwFWFww,
$$
$$
FFwfwFwfWfwfWWfwfWfwFwfwFFwwFWFww>
$$

$$f=
\left(
  \begin{array}{cc}
    1.40427-1.17926 i & 0 \\
    0 & 0.417611+0.350696 i \\
  \end{array}
\right)
$$
$$
w=\left(
  \begin{array}{cc}
    1.07481-0.850372 i & 0.313498-1.03464 i \\
    0.493763-0.322133 i & 0.747073+0.0218061 i \\
  \end{array}
\right)
$$
According to [L] and Snap the volume of $N_3$ is 7.73809... So, the manifold
$N_3$ possibly could be a 2, 3, 4, 5, 6, 7 and 8-fold cover of a non
exceptional manifold. Therefore, a shortest closed geodesic $\kappa$
of the underlying manifold $M$ must have a tube of radius not less
than $\log{3}/2$. Consequently, its preimage must be a closed geodesic
$\delta$ (or several geodesics $\{\delta_i\}$) with tube radii $>\log(3)/2$. The density
of a volume of a tube $W$ of radius $r\geq \log(3)/2$ about a
geodesic $\delta$ must be less than or equal to the $0.91$ by [P] and hence volume$(W)\le 7.04166$. This implies that  any geodesic of tube radius $>\log(3)/2$ must have length at most $ 6.72429$. From Snap we obtain the list of all
geodesics of length up to 6.7243 for the manifold $N_3$.  Again we do not display geodesics with ortholines of real length $<\log(3)$.

\begin{table}[htdp]
\begin{center}
Length spectrum for $N_3$
\begin{tabular}{|c|c|c|c|}
\hline
orbit & geodesic length & shortest ortholine & geodesic number\\
\hline
0 & 1.21275-1.39704*i & 1.09488+1.17345*i & 0,1,2,3*,4*,5* \\
\hline
1 &  1.59139+2.39677*i & \bf{1.67039-2.41832*i} & 6,7,8,9*,10*,11* \\
\hline
2 & 1.94977-2.59941*i & 0.82700+0.48158*i & 12,13,14,15*,16*,17* \\
\hline
3 &  2.59953-0.00000*i & \bf{1.29867-2.03065*i} & 18,21,19,22,23,20 \\
\hline
 ... & & & \\
\hline
\end{tabular}\\
Length of a geodesic marked by asterisk is conjugated to the listed length.
\end{center}
\end{table}

There are two orbits [1] and [3]
of six geodesics  each of real length 1.59139 and 2.59953 and with tube radii respectively
of 0.835193 and 0.649334. The ratio of these lengths is not a rational number $a/b$ such that $a+b=n$, $n=2,3,4,5,6,7$. Therefore the preimage of $\kappa$ must be
either in orbit [1] or [3].

It cannot happen that $\delta_i$ is degree-1 cover of $\kappa$
because $\kappa$ is a shortest geodesic in $M$ but an image of a geodesic
from orbit [0] is shorter than $\delta_i$.

Let $\delta$ denote an element of orbit [1].  Assume that $p:N_3\to M$ and $p|\delta$ is a degree-$x$ cover of $\kappa$ where $\kappa$ is a shortest geodesic in $M$ and $x>1$.

\subsubsection{$\delta\in$ Orbit [1]:}

 Snap asserts that there is a unique ortholine $\sigma$ of real length 1.81586... that starts and ends on $\delta$.  Furthermore, there is no ortholine of that real length between distinct geodesics $\delta$ to $\delta'$ from orbit [1].
This contradicts the fact that preimage of a shortest ortholine of $\kappa$ has $x$ ortholines between $\delta$ and $\delta'$.

{\footnotesize
\begin{table}[htdp]
\begin{center}
Ortholines of a geodesic from orbit [1] for $N_3$
\begin{tabular}{|c|c|c|}
\hline
ortholine length & initial point of ortholine & final point of ortholine\\
\hline
 1.67039-2.41832*i  &  0.22820-2.70454*i  & 0.39785+1.07038*i \\
 \hline
 1.67039-2.41832*i  & -0.39785-0.12800*i  & -0.56749+2.38026*i \\
 \hline
 \bf{1.81586+0.21843*i}  &  0.31302-0.81708*i  & -0.48267-2.01546*i \\
 \hline
 1.85170-0.80441*i  & 0.18786-1.59286*i  & 0.43819-0.04130*i \\
 \hline
 1.85170-0.80441*i &-0.35751-1.23968*i  & -0.60783-2.79125*i  \\
\hline
 ... & & \\
\hline
\end{tabular}
\end{center}
\end{table}%
}


{\footnotesize
\begin{table}[htdp]
\begin{center}
\begin{tabular}{|c|c|}
\hline
ortholine length & multiplicity\\
between $\delta_1$ and $\delta_2$ from [1] & \\
\hline
 1.15910-0.00000*i & 2\\
\hline
 1.15910+3.14159*i & 2 \\
\hline
 2.17052-0.51494*i & 2\\
\hline
 2.17052+0.51494*i & 2 \\
\hline
 2.17052-2.62666*i & 2 \\
\hline
 2.17052+2.62666*i & 2  \\
\hline
 ... & \\
\hline
\end{tabular}
\end{center}
\end{table}%
}

\subsubsection{$\delta\in$ Orbit [3]:}

By Snap any two distinct geodesics from orbit [3]  have an ortholine connecting them of
 real length less than $\log(3)$. This implies that
only a single geodesic
from orbit [3] can be a preimage of $\kappa$, hence we can assume that that geodesic is $\delta$. Observe that $x> 2$  since $\delta$ is more than two
times longer than a shortest geodesic of $N_3$. On the other hand, by Snap
$\delta$ has only two shortest ortholines that begin and end on $\delta$.  This implies that $x\le 2$ a contradiction.

{\footnotesize
\begin{table}[htdp]
\begin{center}
Orthospectrum of orbit [3] for $N_3$
\begin{tabular}{|c|c|c|}
\hline
ortholine length & initial point of ortholine & final point of ortholine\\
\hline
 1.29867-2.03065*i & -0.03493-0.46561*i & -1.26483-2.59858*i \\
\hline
 1.29867+2.03065*i & 0.03493+2.67598*i & 1.26483+0.54301*i \\
\hline
 1.29909-1.04518*i & 0.02157+1.04332*i & 1.27820+2.17567*i \\
\hline
 1.29909+1.04518*i & -0.02157-2.09827*i & -1.27820-0.96592*i \\
\hline
 ... & & \\
\hline
\end{tabular}
\end{center}
\end{table}%
}

It means there is no chance that any closed geodesic or geodesics
from orbits [1] or [3] can map onto a shortest geodesic $\gamma$ of
$M$.

We now consider the case that $M$ is an exceptional manifold. As in the previous section we can assume that no other exceptional manifold can cover $M$.
 Thus a geodesic $\gamma$ from orbit [0] of $N_3$ maps onto a shortest geodesic $\kappa$ of $M$. Hence real length $\kappa>1.21275$. The preimage of a shortest ortholine $\omega$ of $\kappa$ is a shortest ortholine of $\gamma$, $Relength(\omega)>1.09488$. The maximum volume of tube $W$ around $\kappa$ is more than  $1.26053$. By [P] the volume of the tube W has density at most 0.91 in $M$ and hence $vol(M)>1.3851$.  Therefore the $\deg(p) \leq 5$.

There are no geodesics in $N_3$ with real length 2, 3, 4 or 5 times bigger than real length of  a single geodesic from orbit [0].  A preimage of $\kappa$ lays entirely in orbit [0].
Real length of shortest ortholines between sets of geodesics \{0, 1, 2\} and \{3, 4, 5\} is 0.12450É . It is less then the shortest ortholine $\omega$ therefore there is no 3, 4, 5-fold cover.
There might be a 2-fold cover. One component of $p^{-1}(\kappa)$ is a geodesic 0, 1 or 2 and the other component is 3, 4 or 5.  But since the complex parts of geodesics 0, 1, 2 differ from 3, 4, 5 by sign, it follows that any isometry that takes 0 to 5 (or 3 or 4)  is orientation reversing.  Which would imply that the quotient is non orientable.  This is a contradiction for homological reasons [CLLMR].

{\footnotesize
\begin{table}[htdp]
\begin{center}
Ortholine spectrum up to 1.9 between geodesics from orbit [0] for $N_3$
\begin{tabular}{|c|c|c|}
\hline
ortholine length & initial point of ortholine & final point of ortholine\\
\hline
0.12450-1.05142*i & 0: 0.44052+1.51057*i & 1: 0.16585-0.31956*i \\
 0.12450-1.05142*i & 1:-0.44052-2.76263*i & 2:-0.16585+0.40414*i \\
 0.12450+1.05142*i & 3: 0.44052+0.37896*i & 4: 0.16585-2.73745*i \\
 0.12450+1.05142*i & 4:-0.44052-0.29438*i & 5:-0.16585-0.63829*i \\
 0.12450-2.09017*i & 3:-0.16585+2.82203*i & 5: 0.44052-3.08136*i \\
 0.12450+2.09017*i & 0:-0.16585-0.93251*i & 2: 0.44052+2.84722*i \\
 1.09488-1.17345*i & 0: 0.13733+0.28903*i & 0:-0.46904-2.15404*i \\
 1.09488-1.17345*i & 1:-0.13733-1.54110*i & 1: 0.46904+0.90198*i \\
 1.09488-1.17345*i & 2: 0.13733+1.62568*i & 2:-0.46904-0.81740*i \\
 1.09488+1.17345*i & 3: 0.13733+1.60049*i & 3:-0.46904-2.23962*i \\
 1.09488+1.17345*i & 4:-0.13733-1.51591*i & 4: 0.46904+2.32420*i \\
 1.09488+1.17345*i & 5: 0.13733-1.85983*i & 5:-0.46904+0.58325*i \\
 1.57206-1.08896*i & 3:-0.02530+1.98594*i & 5: 0.29996-2.24527*i \\
 1.57206-1.08896*i & 3:-0.30641-2.62506*i & 5: 0.58108+2.36573*i \\
 1.57206+1.08896*i & 0:-0.02530-0.09642*i & 2: 0.29996+2.01113*i \\
 1.57206+1.08896*i & 0:-0.30641-1.76860*i & 2: 0.58108-2.59988*i \\
 1.57206-2.05263*i & 0: 0.29996+0.67448*i & 1: 0.02530-1.15565*i \\
 1.57206-2.05263*i & 0: 0.58108+2.34666*i & 1: 0.30641+0.51653*i \\
 1.57206-2.05263*i & 1:-0.29996-1.92655*i & 2:-0.02530+1.24023*i \\
 1.57206-2.05263*i & 1:-0.58108+2.68446*i & 2:-0.30641-0.43195*i \\
 1.57206+2.05263*i & 3: 0.29996+1.21505*i & 4: 0.02530-1.90136*i \\
 1.57206+2.05263*i & 3: 0.58108-0.45713*i & 4: 0.30641+2.70965*i \\
 1.57206+2.05263*i & 4:-0.29996-1.13047*i & 5:-0.02530-1.47438*i \\
 1.57206+2.05263*i & 4:-0.58108+0.54171*i & 5:-0.30641+0.19780*i \\
 ... & & \\
\hline
\end{tabular}
\end{center}
\end{table}%
}

\subsection{Manifold N4}
For the region $X_4$ there exists a unique manifold $N_4$. The fundamental group of this manifold
is:
$$
<f,w\,|\;FFwfwFwfWfwfWfwFwfwFFwwFWFwFWFww,
$$
$$
FFwfwFwfwFFwwFWFwFWfWFWfWFwFWFww>
$$

$$f=
\left(
  \begin{array}{cc}
    1.35462-1.22513 i & 0 \\
    0 & 0.40607+0.367252 i \\
  \end{array}
\right)
$$
$$
w=\left(
  \begin{array}{cc}
    1.02306-0.877334 i & 0.265945-1.07164 i \\
    0.501555-0.337493 i & 0.737634+0.0194601 i \\
  \end{array}
\right)
$$
We repeat examination of the manifold $N_4$ in the same way as for previous manifolds.
According to [JR] and Snap the volume of $N_4$ is 7.517689... So, the
$N_4$ possibly could be a 2, 3, 4, 5, 6 and 7-fold cover of a non
exceptional manifold $M$. Therefore, a shortest closed geodesic $\kappa$
of the underlying manifold $M$ must have a tube of radius not less
than $\log{(3)}/2$. Consequently, its preimage must be a closed geodesic
$\delta$ (or some geodesics $\delta_i$) with thick tuberadius. The volume of a tube $W$ of radius $r\geq \log(3)/2$ about a
geodesic $\delta$ must be less or equal to the $0.91 Vol(N_4)=6.8411$ by [P]. We see that all thick tube geodesics
must be of length $l\leq 6.53277$. From Snap we get the list of all
geodesics of length up to 6.53277 for manifold $N_4$. Geodesics which do not appear in the table have shortest ortholines less than $\log(3)$.

\begin{table}[htdp]
\begin{center}
Geodesics for $N_4$
\begin{tabular}{|c|c|c|c|}
\hline
orbit & geodesic length & shortest ortholine & geodesic number\\
\hline
  0 & 1.20475+1.47049*i & 1.09508+1.23769*i & 0,1,2,3,4,5 \\
\hline
 1 & 1.36612-2.17271*i & \bf{1.90660+2.73377*i} & 6,7,8 \\
\hline
 2 &  1.90660+2.73377*i & 0.86339+0.51521*i & 9,10,11,15,16,17 \\
 \hline
 3 & 1.90660+2.73377*i & \bf{1.36612-2.17271*i} & 12,13,14 \\
\hline
 4 & 2.57004+0.51982*i & \bf{1.16156-1.35749*i} & 18,19,20,21,22,23 \\
\hline
 ... & & & \\
\hline
 18 & 3.79966-1.28485*i & \bf{1.16280+1.53153*i} & 164,165 \\
\hline
 ... & & & \\
\hline
\end{tabular}
\end{center}
\end{table}

There are four orbits [1], [3], [4] and [18]
of three, three, six and two geodesics of real length 1.36612, 1.90660, 2.57004 and 3.79966 with injectivity radii 0.953299, 0.68306, 0.580779 and 0.581399. Ratios of these lengths do not give rational numbers $a/b$ such that $a+b=n$, where $n=2,3,4,5,6,7$. Therefore the preimage of $\kappa$ must lie
either in orbit [1], [3], [4] or [18].

There is no possibility that $\delta_i$ is a one-to-one cover of $\kappa$ because $\kappa$ is a shortest geodesic but the image of geodesics
from orbit [0] is shorter than $\kappa$.

\subsubsection{Orbit [1]}

Shortest ortholines between different geodesics from orbit [1] have real length $0.95330<\log(3)$. That is the reason why we can consider as a preimage of $\kappa$ only a single geodesic $\delta$
from orbit [1].
Each geodesic $\delta_i$ from orbit [1] has only two shortest ortholines of length
$1.906$. These ortholines cannot be a preimage of any ortholine of
$\kappa$ under the projection $p$ except if $p$ is a double cover and $\delta_i$ is a double cover of $\kappa$.

{\footnotesize
\begin{table}[htdp]
\begin{center}
Ortholine spectrum for a geodesic from orbit [1] for $N_4$
\begin{tabular}{|c|c|c|}
\hline
ortholine length & initial point of ortholine & final point of ortholine\\
\hline
 1.90660+2.73377*i  & 0.18220-0.99465*i  & 0.18220+2.14695*i \\
\hline
 1.90660+2.73377*i  & 0.86526-2.08100*i  & 0.86526+1.06059*i \\
\hline
 2.04203+0.92305*i  & 0.10760+1.27981*i  & 0.25681-0.12751*i \\
\hline
 2.04203+0.92305*i  & 0.10760-1.86178*i  & 0.25681+3.01408*i \\
\hline
 2.04203+0.92305*i  & -0.42625+0.95884*i  & 0.79066+0.19346*i \\
\hline
 2.04203+0.92305*i  & -0.42625-2.18275*i  & 0.79066-2.94813*i  \\
\hline
 ... & & \\
\hline
\end{tabular}
\end{center}
\end{table}%
}
We check all closed
orientable hyperbolic manifolds from SnapPea's census which can be potentially a double
quotient of $N_4$ using volume and the first homology group
$H_1=\mathbb{Z}_4\oplus \mathbb{Z}_{12}$. We get only one such
manifold $m371(1,3)$. One of
its double covers $m371\sim 2(1,0)$ has a fundamental group
$$
<a,b|\,aabbaBabaBabbaabAbabABAbABAbabAb,
$$
$$
aabbaabAbabABAbABBAABBAbABAbabAb>
$$
which is isomorphic to the
$\pi_1(N_4)$:
$$
a\rightarrow w,\quad b \rightarrow F,\quad \mbox{and} \quad w\rightarrow a, \quad f\rightarrow B.
$$
This proves the existence of a double quotient of
$N_4$.

In a similar way as we did for $N_2$ we can show that $N_4$ does not have any other quotients. As for $N_2$ we see from orthospectrum of orbit [1] that there are at most three free $\mathbb{Z}_2$ symmetries of $N_4$.
Geodesics of orbit [1] have the following representations in terms of generators of the fundamental group:
$$6\rightarrow FwwF, \quad 7\rightarrow WFwF, \quad 8\rightarrow WfWF.$$
From Snap we got three symmetries that map an element of the orbit [1] to any other element of the same orbit:
\begin{enumerate}
\item $f\rightarrow WffWWfwfW,\quad w\rightarrow WffWFWfWF,$
\item $f\rightarrow wFWFwFWfWF,\quad w\rightarrow wFWFwwFFw,$
\item $f\rightarrow fwFwfwFFw,\quad w\rightarrow fwFwfWfwfW.$
\end{enumerate}
We rigorously checked this by hand by finding automorphisms of $\pi_1(N_4)$ that have this property. It follows that a projection $q$ is conjugate to $p$ with the conjugating map being one of these three isometries.

{\footnotesize
\begin{table}[htdp]
\begin{center}
Ortholines between different geodesics of orbit [1] for $N_4$
\begin{tabular}{|l|l|l|}
\hline
ortholine length & initial point of ortholine & final point of ortholine\\
\hline
 0.95330+1.36689*i & 6: 0.86526-2.08100*i & 7: 0.47302+1.69236*i  \\
\hline
 0.95330+1.36689*i & 6: 0.86526+1.06059*i & 7: 0.47302-1.44923*i  \\
\hline
 0.95330+1.36689*i & 7:-0.21004-0.36288*i & 8: 0.65655-1.18161*i  \\
\hline
 0.95330+1.36689*i & 7:-0.21004+2.77871*i & 8: 0.65655+1.95999*i  \\
\hline
 0.95330-1.77471*i & 6: 0.18220-0.99465*i & 8:-0.02651+3.04634*i  \\
\hline
 0.95330-1.77471*i & 6: 0.18220+2.14695*i & 8:-0.02651-0.09525*i \\
\hline
 ... & & \\
\hline
\end{tabular}
\end{center}
\end{table}%
}

\subsubsection{Orbit [3]}

The distance between distinct geodesics from orbit [3] is $0.68306<\log(3)$ and a preimage of $\kappa$ can be only a single geodesic $\delta$
from orbit [3].  The only possible quotients of $\delta_i$ that must be
considered are 3, 4, 5, 6 and 7. There is no 2-fold cover since $\delta$ is more than two
times longer than a shortest geodesic of $N_4$. Each geodesic
$\delta$ from [3] has only two shortest ortholines of length $1.36612$. Therefore these ortholines cannot be
preimages of any underlying ortholine of $\kappa$.

{\footnotesize
\begin{table}[htdp]
\begin{center}
Ortholine spectrum for a geodesic from orbit [3] for $N_4$
\begin{tabular}{|c|c|c|c|}
\hline
 & ortholine length & initial point of ortholine & final point of ortholine\\
\hline
0 & 1.36612-2.17271*i & 0.43080+1.52064*i &  0.43080-1.62096*i \\
\hline
1 & 1.36612-2.17271*i & 1.38410+2.88752*i & 1.38410-0.25407*i \\
\hline
2 & 1.93611+2.16313*i & 0.33808-0.37206*i & 1.29137+0.99483*i \\
\hline
3 & 1.93611+2.16313*i & 0.33808+2.76954*i & 1.29137-2.14676*i \\
\hline
4 & 1.93611+2.16313*i & -0.42977-1.09515*i & 0.52353+0.27174*i \\
\hline
5 & 1.93611+2.16313*i & -0.42977+2.04644*i & 0.52353-2.86986*i \\
\hline
 ... & & &\\
\hline
\end{tabular}
\end{center}
\end{table}%
}

\subsubsection{Orbit [4]}

The distance between distinct geodesics from orbit [4] is $0.15119$, $0.57139$ or $0.68605$, all less than $\log(3)$. A preimage of $\kappa$ can be only a single geodesic $\delta$
from orbit [4].  Each geodesic $\delta_i$ from orbit [4] has two shortest ortholines of length
$1.16156$. Then we can apply the same line of argument as we did before.

{\footnotesize
\begin{table}[htdp]
\begin{center}
Ortholine spectrum for a geodesic from orbit [4] for $N_4$
\begin{tabular}{|c|c|c|c|}
\hline
 & ortholine length & initial point of ortholine & final point of ortholine\\
\hline
0 & 1.16156-1.35749*i & -0.15089-0.15348*i & 0.86832+1.39935*i \\
\hline
1 & 1.16156-1.35749*i & -0.41671-2.00215*i & 1.13413-3.03516*i \\
\hline
2 & 1.44803+1.88714*i & -0.11882+2.72203*i & 0.83625-1.47616*i \\
\hline
3 & 1.44803+1.88714*i & -0.44877+1.40552*i & 1.16620-0.15965*i \\
\hline
 ... & & &\\
\hline
\end{tabular}
\end{center}
\end{table}%
}

\subsubsection{Orbit [18]}

 The distance between different geodesics from orbit [18] is $0.70700<\log(3)$ and a preimage of $\kappa$ can be only a single geodesic $\delta$
from orbit [18].
Possible quotients of $\delta_i$ are 4, 5, 6 and 7. There are no 2 and 3-fold covers since $\delta$ is more than three times longer than a shortest geodesic of $N_4$. Each geodesic
$\delta$ from [18] has only three ortholines of length 1.23288. Hence these ortholines cannot be
preimages of any underlying ortholine of $\kappa$.

{\footnotesize
\begin{table}[htdp]
\begin{center}
Ortholine spectrum up to 1.78 for a geodesic from orbit [18] for $N_4$
\begin{tabular}{|c|c|c|}
\hline
 ortholine length & initial point of ortholine & final point of ortholine\\
\hline
1.16280+1.53153*i & -0.21988+1.05297*i & 0.30723-2.18999*i \\
\hline
 1.16280+1.53153*i & -0.21988-1.04143*i & 0.30723+1.99880*i \\
\hline
 1.16280+1.53153*i &-0.21988-3.13582*i & 0.30723-0.09559*i \\
\hline
 1.16280+1.53153*i &-1.59260+0.54683*i & 1.67995+2.50494*i \\
\hline
 1.16280+1.53153*i & -1.59260-1.54756*i & 1.67995+0.41054*i \\
\hline
 1.16280+1.53153*i & -1.59260+2.64123*i & 1.67995-1.68385*i \\
\hline
 1.23288-1.97709*i & -0.90624-0.24730*i & 0.99359-0.88972*i \\
\hline
 1.23288-1.97709*i & -0.90624+1.84710*i & 0.99359+1.20467*i \\
\hline
 1.23288-1.97709*i & -0.90624-2.34169*i & 0.99359-2.98412*i \\
\hline
 ... & & \\
\hline
\end{tabular}
\end{center}
\end{table}%
}
It means there is no chance that a closed geodesic or geodesics
from orbits [1], [3], [4] or [18] can map onto a shortest geodesic $\kappa$ of
$M$ except if $M$ is the manifold $m371(1,3)$.

Finally we consider the case when $M$ is an exceptional manifold.
A geodesic $\gamma$ from orbit [0] of $N_4$ maps onto a shortest geodesic $\kappa$ of $M$ and real length $\kappa>1.20475$. The preimage of a shortest ortholine $\omega$ of $\kappa$ is a shortest ortholine of $\gamma$. Thus $\Relength(\omega)>1.09508$. The maximum volume of tube $W$ around $\kappa$ is more than  $1.25272$. By [P] the volume of the tube W has density at most 0.91 in $M$ and hence $vol(M)>1.37661$.  Therefore the $\deg(p) \leq 5$.

There are no geodesics of length 2, 4 or 5 times of $\gamma$. There are two orbits [14], [15] of geodesics  with real length three times bigger than length of  a single geodesic $\gamma$ from orbit [0]. Their shortest ortholines are of length 0.90630É and 1.07165É . Both are less than $\Relength \omega$. Therefore a preimage of $\kappa$ lays entirely in orbit [0].
Real length of shortest ortholines between set of geodesics \{0, 1, 3\} and \{6, 7, 8\} is 0.13192. It is less then the shortest ortholine $\omega$ therefore there is no 3-fold cover.
There might be a 2-fold cover. One component of $p^{-1}(\kappa)$ is a geodesic 0, 1 or 3 and the other component is 6, 7 or 8.  Orbit [1] has three geodesics of real length $1.36612É$.  At least one of them will cover a geodesic of real length 0.68306 that is less than the length of shortest geodesic $1.20475$. Therefore, $N_4$ nontrivially cover no exceptional manifold if the length and ortholength spectra are rigorous.

{\footnotesize
\begin{table}[htdp]
\begin{center}
Ortholine spectrum up to 1.57 between geodesics from orbit [0] for $N_4$
\begin{tabular}{|c|c|c|}
\hline
ortholine length & initial point of ortholine & final point of ortholine\\
\hline
 0.13192+1.05846*i & 0: 0.53079-1.42015*i & 1: 0.24072+1.64443*i \\
 0.13192+1.05846*i & 1:-0.36165-2.23240*i & 2:-0.56763+0.84640*i \\
 0.13192+1.05846*i & 3: 0.20134+1.55778*i & 4: 0.52205-2.98266*i \\
 0.13192+1.05846*i & 4:-0.08032-0.57631*i & 5: 0.66320+2.76679*i \\
 0.13192-2.08313*i & 0:-0.07159+0.98620*i & 2: 0.03475-1.55995*i \\
 0.13192-2.08313*i & 3: 0.80372-0.84857*i & 5: 0.06082-1.11004*i \\
 1.09508+1.23769*i & 0: 0.22960-0.21698*i & 0:-0.37278+2.18937*i \\
 1.09508+1.23769*i & 1:-0.06047+2.84761*i & 1: 0.54191+0.44126*i \\
 1.09508+1.23769*i & 2:-0.26644-0.35677*i & 2: 0.33593-2.76312*i \\
 1.09508+1.23769*i & 3:-0.09985+2.76096*i & 3: 0.50253+0.35461*i \\
 1.09508+1.23769*i & 4: 0.22086-1.77949*i & 4: 0.82324+2.09735*i \\
 1.09508+1.23769*i & 5:-0.24037+0.09313*i & 5: 0.36201-2.31322*i \\
  ... & & \\
\hline
\end{tabular}
\end{center}
\end{table}%
}

\section{Rigorous length and ortholength spectra}

In the previous section we used results obtained by Snap: Dirichlet domain, length and
ortholength spectra, injectivity radii to prove the theorem. That part was based on the experimental data which are not rigorous (the round-off error was not considered). We took exactly the same input from Snap (face parings) and found the list of geodesics and ortholines with exact round-off errors (we did all algebraic calculations with errors). Our results are identical to ones from Snap that proves the theorem.

To check the results from Snap we write a package in {\sc Mathematica} that calculates the length and ortholength spectra for the given geodesic. In this section we describe the theoretical part of our algorithm. Two files are attached to the arxiv version of this paper - {\it length.nb} (which contains the interface where you actually run the code together with the description of all necessary commands and options), {\it source$\_$length.nb}  includes a source code which is loaded by the first file in the beginning. Our package was written in order to find the geodesics and ortholines for manifolds $N_2$, $N_3$, $N_4$ but can be used for any other manifolds as well.

As for precision, Mathematica uses advanced algorithms to reach arbitrary precision during numerical evaluations. Therefore, the precision of our result is limited just by the precision of the input data (generators of a fundamental group, face parings of a Dirichlet domain) and the computer memory. Also, the precision of the evaluation is being continuously updated when running the code. Our algorithm was tested on the cases of manifolds $N_2$, $N_3$, $N_4$ (and also other examples) and the results for geodesics and ortholines precisely agree with the data obtained by Snap version 1.11.3. For each manifold we
calculate geodesic length, injectivity radius, etc. These numbers
are never precise because the initial data have some error. But
because we keep track of the errors at all stages we know the
round-off error also for these numbers. Therefore, it is trivial to
check if the given number is smaller or bigger than some cutoff
value (called precision in the code) within its precision (round-off error). Therefore, we can
claim that our result is rigorous.

For manifold $N_3$ we want to calculate geodesics up to length 6.7243 and their tube radii. Dirichlet domain for $N_3$ has 36 faces. A lot of computations must be done to get the list of geodesics in this case. The number of order $10^8$ group elements must be checked before we get a list of geodesics. It takes several hours for package to get the list while Snap calculates this list in few minutes, but reliability is our first priority. For the same manifold the package gives a list of geodesics with cut-off less than 4.0 just in seconds.

We remind that manifold $N_i$ is a unique exceptional manifold associated to the region $X_i$ and manifold $M_i$ is calculated by Snap using fundamental group of $N_i$.
There is an isomorphism between fundamental
groups of manifolds $N_i$ and $M_i$, $i=\{2,3,4\}$. From Mostow's Rigidity Theorem [BP] it follows that these manifolds are homeomorphic.

For example, the fundamental group of $M_2$ is
$$
< a b |\, AbabaBaaBababAbb,\; AAbAAbbAbababAbb >
$$
and an isomorphism is $a\rightarrow f, b\rightarrow w$.

Let $M$ be a hyperbolic three-dimensional manifold of finite volume. We use information about a Dirichlet domain of $M$ from Snap to calculate the
spectra and injectivity radii precisely. We use:

\begin{itemize}
\item  fundamental group represented by generators $a_i$ as matrices $SL(2, C)$ with high precision,
\item a matrix $c\in O(3,1)$ conjugating between coordinates in which the base point of the Dirichlet domain is at the origin and the coordinates in which the generators for the fundamental group were originally given to the Dirichlet domain finding code,
\item words for the face parings.
\end{itemize}

An algorithm for length spectrum is described in [HW]. The difference of our algorithm from the one used in Snap is that we construct the Dirichlet domain for $M$ in the projective ball
model (also known as Beltrami's model or Klein's model). In this model we operate with matrices $SL(2,C)$ which give smaller error than $O(3,1)$ and it is easier to calculate hyperbolic distance, edges and planes there.

The algorithm requires a spine radius of Dirichlet domain. The \textit{spine radius} is defined as the infimum of the radii of all spines to the domain. We cannot get a spine radius of the Dirichlet domain from Snap. Therefore, we start with a construction of the Dirichlet domain for $M$ and calculating a lower bound of its spine radius.

We recall the definition of the Dirichlet domain [Bo]. Let $\Gamma$ be a group of isometries of the metric space $(X,d)$ whose action is discontinuous. The \textit{Dirichlet domain} of $\Gamma$ centered at the point $x\in X$ is the subset
$$
D_\Gamma (x)=\{y\in X: d(x,y)\leq d(g (x),y) \mbox{ for every } g\in\Gamma\},
$$
consisting of those points $y\in X$ which are at least as close to $x$ as to any other point of its orbit $\Gamma (x)$. In our notation we refer to $D_\Gamma (x)$ as $D$ and always $x$ will be a basepoint. In the case of hyperbolic space $H^3$ the set of points $z$ that are at the same hyperbolic distance from $x$ and $y$ is a hyperbolic plane $P_g$ and the set of $z$ with $d_{{\footnotesize \mbox{hyp}}}(x,z)\leq d_{{\footnotesize \mbox{hyp}}}(y,z)$ is a hyperbolic half-space $H_g$ delimited by this perpendicular bisector plane $P_g$. The Dirichlet domain $D$ is a finite-sided polyhedron. The polyhedra $gD$ with $g\in\Gamma$ form a tessellation of $H^3$ and $gD$ is distinct from $D$ unless $gx\ne x$.

For our purposes it is better to define the Dirichlet domain using half-spaces. The Dirichlet domain $D$ of $M$ with base point $x$ is the intersection of the half-spaces $H_g$, for all covering transformations $g\in \Gamma$:
$$
D=\bigcap_{g\in \Gamma} H_g.
$$

Each hyperbolic isometry $g$ has an axis $A_g$, that is a fixed geodesic under the isometry. In other words, each transformation $g$ corresponds a geodesic $A_g$. If we look for geodesics of length up to $\lambda$ we have to consider all isometries $g$ that move the basepoint $x$ a distance less than $s$. This distance $s$ depends on the cut-off length $\lambda$ and the size of the Dirichlet domain, which is characterized by the spine radius.

\subsection{Spine Radius}

After giving some definitions we define a spine radius and introduce an algorithm how to calculate it.

Each Dirichlet domain, with faces identified, specifies a cell decomposition $K$ for $M$. A \textit{spine dual to the Dirichlet domain} is a two-skeleton of a cell decomposition $K'$ of $M$ dual to $K$. All closed geodesics of $M$ intersect a spine dual to the Dirichlet domain. The maximum distance from a point in the spine dual to the Dirichlet domain called its \textit{radius}.

\textit{The spine radius $r$ of the Dirichlet domain} is the infimum of the radii of all spines dual to the domain.

The spine radius is the maximin edge distance of the Dirichlet domain [HW]. We use this fact to calculate a spine radius rather than its direct definition. Note that the spine radius is finite for all Dirichlet domains.

We explain the algorithm of the construction of the Dirichlet domain for a manifold $M$ and use data about the Dirichlet domain obtained by Snap.

\begin{itemize}
  \item First, we transform $SO(3,1)$ matrix $c$ into $SL(2,C)$  
  and conjugate matrices $a_i$ by the matrix $c$ but continue to refer to them as $a_i$. After conjugating the Dirichlet domain is centered at the point $O(0,0,1)$ in the upper-half space model $U^3$.
\\
\item Represent all face paring relations as matrices $g_j\in SL(2,C)$ (e.g., $j=\overline{1,24}$ for the manifold $M_2$).
\\
\item Find the images $g_j(O)$ of the basepoint $O(0,0,1)$ under all $g_j$ in the $U^3$. They correspond to the basepoints of all neighbor domains. We calculate $g_j(O)$ as a multiplication of quaternions:
\begin{equation}
g(w) = (\alpha\star w + \beta)\star (\gamma\star w+\delta)^{-1},
\nonumber
\end{equation}
where $g=
\left(
\begin{array}{cc}
 \alpha & \beta  \\
 \gamma  & \delta
\end{array}
\right)
\in SL(2,C),
$
$w=x+yi+zj$ and $\star$ represents a multiplication of quaternions.
\\
\item In the paper of [L] it is described how to get the set of all vertices of the Dirichlet domain. We follow this algorithm. Map $O$ and $g_j(O)$ from the upper-half space $U^3$ into projective ball model $D^3$. There is an isometry $\varphi: D^3\longrightarrow U^3$ with inverse $\varphi^{-1}: U^3\longrightarrow D^3$:
$$
\varphi(x+y\,i+z\,j)=\frac{(2x_0,\;2y_0,\;x_0^2+y_0^2+z_0^2-1)}{1+x_0^2+y_0^2+z_0^2}
$$
\\
\item All planes $P_g$ which contain faces of the Dirichlet domain are bisecting planes between points $(0,0,0)$ and $\varphi^{-1}g_j(O)$.
$$
P_g=\{r\in D^3: \mathbf n\cdot r=t\}
$$
where $\mathbf n=\varphi^{-1}g_j(O)$, $t=1-\sqrt{1-|\varphi^{-1}g_j(O)|^2}$.
\\
\item A vertex $R=(r_1,r_2,r_3)$ defined by the intersection of three planes as a simultaneous solution of $\mathbf n_1\cdot r=t_1$, $\mathbf n_2\cdot r=t_2$ and $\mathbf n_3\cdot r=t_3$. We discard vertices which lie "above" all planes $P_g$.
\\
\item Each edge defined by a couple of vertices. We check all pairs of vertices if they define an edge of the Dirichlet domain. We discard a line determined by two vertices if it lies in a single plane.
\\
\item For each edge we calculate the distance from the base point (the origin) to the edge. If the distance to the line containing the edge is smaller than a distance to its vertex then the distance to the closest vertex will be the distance to the edge.
Finally, define the spine radius $r$ as the maximum of these distances.
\end{itemize}

\subsection{Geodesics}
Now we are ready to sketch out the algorithm of computing a length spectrum of geodesics described in [HW]. The general idea of our algorithm is the same but some steps were made differently. We overview the algorithm here with some additions.

\medskip

\textbf{Proposition} (Hodgson, Weeks) \textit{To find all closed geodesics of length at most $\lambda$, it suffices to find all translates $gD$ such that $d(x,gx)\le 2\,\cosh^{-1}(\cosh r \cosh (\lambda/2))$.}

Here $r$ is the spine radius for the Dirichlet Domain $D$ centered at the point $x$.

\medskip
The metric $d$ on the upper-half space is given by
$$
\cosh d(x,y)=1+\frac{|x-y|^2}{2x_3y_3}.
$$
We tile a region in $H^3$ around the Dirichlet domain $D$ centered at the origin by all translates $gD$. These translations move the basepoint to a distance less than the distance defined in the above proposition.
We are not interested in group elements $g$ with $Relength(g)=0$ or $Relength(g)>\lambda$ or whose axis do not pass within a distance $r$ of the basepoint (every geodesic must intersect a spine radius $r$). The $Relength$ of a transformation $g$ is the real part of the complex number
$$
\length (g)=2\Arccosh \frac{tr (g)}{2}.
$$
The distance $r$ from the basepoint to the axis of the isometry $g$ is
$$
r=\Arccosh\sqrt{\frac{\cosh d-\cos t}{\cosh s-\cos t}},
$$
here $d=g(0,0)$, $s+i t=\length(g)$.

We want to find all geodesics that satisfy three constraints described above. The main idea of our algorithm is following

\begin{itemize}
\item All the geodesics up to the cut-off $\lambda$ correspond to group elements which can be constructed from the products of the face pairing relations of the domain. This multiplication forms a natural tree-like structure.
\item We move in this tree and for any point check if it passes the base point distance constraint. If the answer is positive, we continue deeper in the tree from this point. Otherwise, we go back.
\item If the geodesic moves the base point to a distance less than $d(x,gx)$ and in addition also other two constraints are satisfied, we check if it is already in the list. If not, we add it.
\item This algorithm is guaranteed to finish in finite time because there is just finite number of group elements of the tree (which correspond to geodesics up to length $\lambda$) that move the base point to a distance smaller than cutoff.
\end{itemize}

There is no element that move the basepoint of the Dirichlet domain a distance less than $s$, all of whose neighbors move the basepoint to a distance greater than $s$ [HW]. Hence our algorithm cannot miss any translation.

This tiling produce a so called \textit{big list} of group elements $g_i$ corresponding to all geodesics of length at most $\lambda$. This list might contain different group elements which correspond to the same geodesic. We want to have precisely one group element in each conjugacy class. Remove group elements that are just powers of others. Discard all conjugates, the inverse and its conjugates for each geodesic. The conjugacy is realized by an element $h$ from the big list such that
$$
d(x,gx)\le 2\,\cosh^{-1}(\cosh r \cosh (\lambda/4)).
$$
We call a \textit{small list} the part of the big list which is left after eliminations of all duplicates. The small list has a length spectrum with correct multiplicities.

\subsection{Ortholines}

We want to find ortholines between closed geodesics $A_f$ and $A_g$ up to length $\delta$ and positions of their endpoints with angles on the geodesics. We look for them among ortholines between preimages of $A_f$, $A_g$ and conjugates to $A_g$ in the universal cover $U^3$.

\begin{enumerate}

 \item
 For easier calculations we map geodesic $A_f$ (- axis of transformation $f$) onto geodesic $B_{0,\infty}$:
 $$
 q:\,A_f \longrightarrow  B_{0,\infty},\qquad q\in \Isom(H^3)
 $$
$B_{0,\infty}$ is an oriented geodesic $\{(0,0,z):\,0<z<\infty \}$ and $(0,0,0)$ is its negative endpoint. The axis of transformation $f'=q^{-1}\,f\,q$ is $B_{0,\infty}$.
 We choose one of transformations $q$:
 \begin{eqnarray}
 q &=&
\left(
\begin{array}{ccc}
  z_1& z_0/(z_1-z_0)  \\
  1 &  1/(z_1-z_0)
\end{array}
\right)\nonumber
\end{eqnarray}
where $\,z_i=\frac{ f_1-f_4\pm\sqrt{(f_1+f_4)^2-4}}{2\,f_3 }$ are endpoints of $A_f$ (i.e. fixed points of $f$ lie on the boundary of $H^3$).
This transformation map geodesic $A_g$ onto axis of transformation $g_0\,=\,q^{-1}\,g\,q$, where $g_0 =
\left(
\begin{array}{cc}
  g0_1& g0_2\\
  g0_3 &  g0_4
\end{array}
\right)$.

Special cases: If $f_3=0$ and $f_1\ne 1$ then
 \begin{eqnarray}
 q &=&
\left(
\begin{array}{ccc}
  1& \frac{f_1f_2}{1-(f_1)^2}  \\
  0 &  1
\end{array}
\right)\nonumber
\end{eqnarray}
If $f_1+f_4=\pm2$ then $f$ is either a parabolic or a pure reflection and does not have an axis. In this case the algorithm will stop.
\\

 \item
 We choose a transformation $h\in \Isom(H^3)$ such that $h:B_{0,\infty}\longrightarrow A_{g_0}$.
 \begin{equation}
 h =
\left(
\begin{array}{ccc}
 z_1& z_0/(z_1-z_0)  \\
  1 &  1/(z_1-z_0)
\end{array}
\right)\nonumber
\end{equation}
where $z_i$ are endpoints of $A_{g_0}$.

Again we have special cases: If $g0_3=0$ and $g0_1\ne 1$ then
 \begin{eqnarray}
 h &=&
\left(
\begin{array}{ccc}
  1& \frac{g0_1g0_2}{1-(g0_1)^2}  \\
  0 &  1
\end{array}
\right)\nonumber
\end{eqnarray}
If $g0_1+g0_4=\pm2$ then $g0$ is either a parabolic or a pure reflection and does not have an axis. In this case the algorithm will stop too.
\\

 \item
 We can calculate distance between geodesics $f$ and $g_0$. It is equal to the complex distance of an orthocurve from geodesic $B_{0,\infty}$ to $A_{g_0}$ because isometries preserve distance. An ortholine between $B_{0,\infty}$ and $A_{g_0}$  is an axis of transformation $k_0=h\tau h^{-1}\tau$ and its length is twice bigger than the distance $d_0$ ($\tau$ is a rotation around geodesic $B_{0,\infty}$). Hence, the distance $d_0$ is defined from the formula [Be, HW, GMT]
$$
\cosh d_0=\cosh \dist(B, h(B))=\otr(h),
$$
where $\otr(h)=h_1h_4+h_2h_3$.
\begin{eqnarray}
\tau &=&
\left(
\begin{array}{ccc}
  i &  0   \\
  0 &  -i
\end{array}
\right) \nonumber
\end{eqnarray}

\item
In a similar way we calculate a distance between $A_f$ and axes of conjugacy of $g$. We conjugate element $g_0$ by group elements $p_i\in \Isom(H^3)$ from the second "big list" and apply transformation $q^{-1}$. These transformations have the same axes as transformations $g_i=(q^{-1}\,p\,q)\,g_0\,(q^{-1}\,p\,q)^{-1}$ which is also the image of the axis $B_0$ under transformations $ph_i$, where $ph_i=q^{-1}p_i\,q\,h$. Hence, distances $d_i$ between $A_f$ and axes of conjugacy of $g$ defined from the formula:
$$
\cosh d_i=\cosh\dist(B, ph_i(B))=\otr(ph_i), \quad (i\geq 1).
$$
Ortholines for these geodesics are axes of transformations $k_i=ph_i\,\tau \,(ph_i)^{-1}\,\tau$.
The following lemma about conjugated group elements says how to get the second "big list".

 \textbf{Lemma.} \textit{If $g_1$ and $g_2$ are two conjugate group elements such that the axis $A_{g_1}$ corresponds to a geodesic within a distance $r$ from the basepoint and axis $A_{g_2}$ is within a distance $\delta$ from a fixed geodesic $f$, then there is a group element $h$ such that $g_2=h\,g_1\,h^{-1}$ and
$$
d (x,hx)\leq \frac{1}{2}(\lambda_f+\lambda_g)+\delta+r.
$$
}

\textit{Proof of Lemma.} We can consider only situation when the
axis $A_f$ coincides with geodesic $B_{0,\infty}$. We are looking
for ortholines that are at distance less than $\lambda_f/2$ from the
basepoint $x$.  Let $Q$ be the end of the perpendicular from the
basepoint $x$ to axis $A_{g_1}$. There are infinitely many covering
transformations that takes $A_{g_1}$ to $A_{g_2}$. We take one $h$
that minimizes the distance between $hQ$ and ortholine's endpoint
$N$. The length of geodesics are $\lambda_1$, $\lambda_2$ therefore
the distance between $N$ and $hQ$ is less or equal  $\lambda_1$.
Then $d (Q,hQ)\leq d (x,M)+d(M,N)+d (N,hQ)+d(hQ,hx)$
$$
d (Q,hQ)\leq=\lambda_f/2+\delta+\lambda_{g_2}+\delta+r.
$$

\begin{figure}
\includegraphics[scale=0.42]{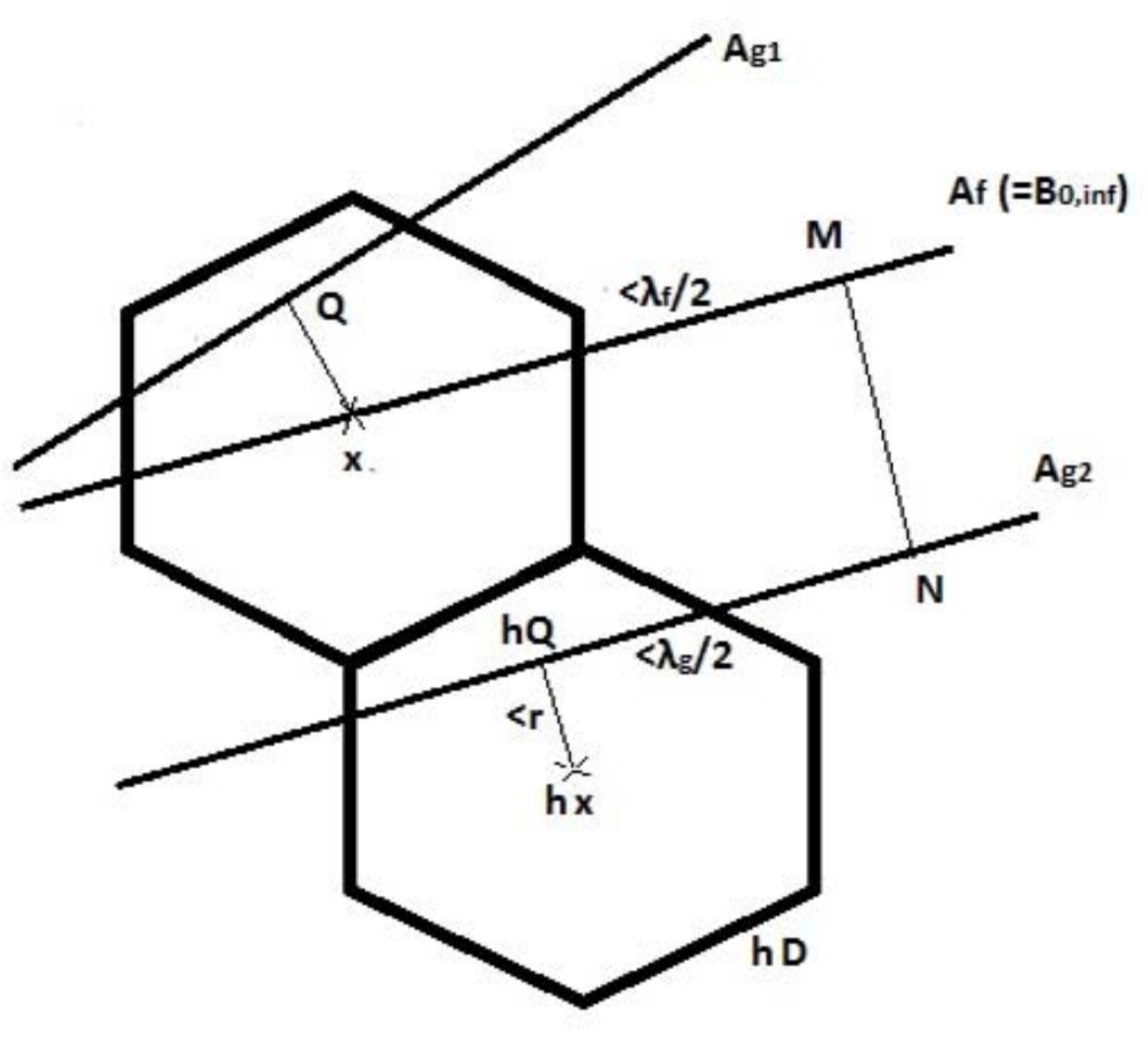}
\caption{The distance from $x$ to $hx$ is
$\leq\frac12(\lambda_f+\lambda_g)+\delta +r$ where $\lambda_\gamma$
is the translation length of $\gamma$.}
\end{figure}

\item
The position of ortholine's endpoint on geodesics $f$ is a complex distance along geodesic $B_{0,\infty}$ from $B_{-1,1}$ to the ortholine between  $B_{0,\infty}$ and $g_j$.
In order to get this distance we fix orientation, choose a point and a based vector on each of these geodesics. The fixed point on the geodesic $B_{0,\infty}$ is the point $(0,0,1)$ of intersection $B_{0,\infty}$ with geodesic $B_{-1,1}$ which runs from endpoint $(-1,0,0)$ to endpoint $(1,0,0)$. The based vector is a tangent vector to the geodesic $B_{-1,1}$ at the fixed point in the positive direction. The positive orientation of geodesic $B_{0,\infty}$ is defined in the direction from endpoint $z_0=(0,0)$ to endpoint $z_1=\infty$ on the boundary $C^2\cup \infty$. This orientation is inherited for all geodesics $g_j$ by covering transformations $ph_j$ (we put $ph_0$ to be $h$ map).  Fixed points on $g_j$ are points of intersection $g_j$ with images of geodesic $B_{-1,1}$ under covering transformations $ph_j$.

Transformations which correspond to ortholines $k_j$ take oriented closed geodesic $f'$ onto oriented closed geodesics $g_j$ along $k_j$ in the positive direction of $k_j$. Endpoints $-T_j$, $T_j$ of ortholines $k_j$ are antipodal to the origin because $k_j$ are perpendicular to geodesic $B_{0,\infty}$. We want to order the endpoints such that ortholines has positive direction at $T_j$. Endpoints $Z_{j0}$, $Z_{j1}$, $-T_j$ and $T_j$ of geodesics $g_j$ and $k_j$ lie on a circle because the geodesics intersect each other. We choose the point $T_j$ to be an endpoint such that point $X=[Z_{j0};\,Z_{j1}]\cap [-T_j;\,T_j]$ belongs to the interval $[0;\,T_j]$. By straightforward calculations this condition is equivalent to the inequality
$$
\left|
\begin{array}{ccc}
 x_0 & y_0  \\
 x_1-x_0 & y_1-y_0
\end{array}
\right|
\left|
\begin{array}{ccc}
 s & t  \\
 x_1-x_0 & y_1-y_0
\end{array}
\right|
\; >\;0.
$$
Here $Z_{j0}=(x_0,\,y_0)=(ph_j)(0)$, $Z_{j1}=(x_1,\,y_1)=(ph_j)(\infty)$, $T_j=(s,t)$.
\\
\item
Now we are ready to calculate the position of the endpoint on the geodesic $f$. The distance from the fixed point on $B_{0,\infty}$ to the orthocurve's endpoint on it is equal to the length of transformation $t_j$ which takes oriented geodesic $B_{-1,1}$ to the oriented orthocurve $k_j$.
Because the ortholine $k_j$ intersects $B_{0,\infty}$ orthogonally then the distance is simply
$$
\length(t_j)=\Log{T_j}
$$
The real part of this formulae $\log |T_j|$ calculates hyperbolic distance between fixed point and the orthocurve's endpoint. The imaginary part equals to the angle between the based vector and tangent vector to the orthocurve at the point $(0,\,0,\,|T_j|)$.

 We consider geodesic $f$ as a circle with period $\lambda=\Relength (f)$. Then ortholine's endpoints on $f$ will be maped onto torus $(-\lambda/2;\,\lambda/2]\times (-\pi i;\,\pi i]$.
\\
\item
Calculations of an endpoint $s_j$ of ortholine $k_j$ on geodesics $g_j$ will be done in a similar way. Transformation $(ph)_j^{-1}$ takes oriented geodesic $g_j$ onto $B_{0,\infty}$ and oriented ortholine $k_j$ onto geodesic which is orthogonal to $B_{0,\infty}$ and have endpoints $(ph)_j^{-1}(T_j)$, $-(ph)_j^{-1}(T_j)$. Then the position is defined by the formulae
$$
\length(s_j)=\Log{(ph)_j^{-1}(T_j)}+\pi i.
$$
We add $\pi i$ because we need an angle between the fixed vector and the orthocurve but not its tangent vector at the endpoint. For geodesics $g_j$, the period is $l=\length (g)$ and endpoints belong to the torus $(-l/2;\,l/2]\times [-\pi i;\,\pi i]$.

\item
We sort all ortholines by three parameters:
\begin{enumerate}
\item distance between geodesics,
\item position of endpoints on geodesics $f$ and $g$.
\end{enumerate}
Then we eliminate duplicate with identical parameters. Geodesics which start at the same point and go in the same direction coincide. We end up with the list of orthocurves up to real length $\delta$ with correct multiplicities.

\end{enumerate}

\textbf{References}

\begin{description}
\item[A] I.~Agol, \textit{Volume change under drilling}, Geom. Top. \textbf{6} (2002), 905-916.
\item[ACS] I.~Agol, M.~Culler, P.~Shalen, \textit{Dehn surgery, homology and hyperbolic volume}
Algebr. Geom. Topol. \textit{6} (2006), 2297Ð2312.
\item[AST] I.~Agol, P.~Storm, W.~Thurston, \textit{Lower bounds on volumes of hyperbolic Haken 3-manifolds}, J. AMS \textbf{20} (2007), no.4, 1053-1077.
  \item [Be] A.~Bearden, \textit{The Geometry of Discrete Groups}, Springer, New York, 1983.
  \item [BP] R.~Benedetti, C.~Petronio, \textit{Lectures in Hyperbolic Geometry}, Springer - Verlag, 1992.
  \item [Bo] F.~Bonahon, \textit{Low-Dimensional Geometry: From Euclidean surfaces to Hyperbolic Knots}, Student mathematical library; v.49. IAS/PARK City mathematical subseries, 2009.
  \item[CLLMR] A.~Champanerkar, J.~Lewis, M.~Lypyanskiy, S.~Meltzer,
  A.~W.~Reid, \textit{Exceptional regions and associated exceptional
  hyperbolic 3-manifolds}, Experiment. Math. \textbf{16} (2007), no. 1, 107- 118.
  \item [CD] M.~Culler, N.~Dunfield, \textit{SnapPy}, Available online from
  http://snappy. computop.org (2009).
  \item[G] D.~Gabai, \textit{On the geometric and topological rigidity of
  hyperbolic 3-manifolds}, J. Amer. Math. Soc. \textbf{10} (1997),
  37-74.
  \item[G1] D.~Gabai, \textit{The Smale conjecture for hyperbolic 3-manifolds: Isom($M^3$)?Diff($M^3$)}. J. Differential Geom. \textbf{58} (2001), no. 1, 113Ð149.
  \item[GM] F.~Gehring, G.~Martin, \textit{Precisely invariant collars and the volume of hyperbolic 3-folds}, J. Differential
Geom. \textbf{49} (1998), 411Ð435.
  \item [GMM] D.~Gabai, R.~Meyerhoff, P.~Milley, \textit{Minimum volume cusped hyperbolic three-manifolds}, J. Amer. Math. Soc. \textbf{22} (2009), no. 4, 1157-1215.
  \item[GMT] D.~Gabai, R.~Meyerhoff, N.~Thurston, \textit{Homotopy
  hyperbolic 3-manifolds are hyperbolc}, Ann. of Math. (2)
  \textbf{157} (2003), no. 2, 335-431.
  \item[GHN] O.~Goodman, C.~Hodgson, W.~Neumann, \textit{Home Page For Snap}, Available online from http://www.ms.unimelb.edu.au/~snap/ (1998).
 \item [HW] C.~Hodgson, J.~Weeks, \textit{Symmetries, isometries and length spectra of closed hyperbolic three-manifolds}, Experiment. Math. Vol.3 (1994), No. 4.
  \item[JR] K.~N.~Jones, A.~W.~Reid, \textit{Vol3 and other exceptional hyperbolic
  3-manifolds}, Proc. Amer. Math. Soc. \textbf{129} (2001), no.
  7, 2175-2185.
 \item [L] M.~Lypianskiy, \textit{A computer-assisted application of poincare's fundamental polyhedron theorem}.
  \item[P] A.~Preworski, \textit{A universal upper bound on density of tube palings in hyperbolic space}, J. Diff. Geom. \textbf{72} (2006), 113-127.
  \item[W] Jeffery R. Weeks, \textit{SnapPea}, Available online from
  http://geometrygames.org/ SnapPea/index.html (1993).
\end{description}

\bigskip
\small
\textsc{Department of Mathematics, Princeton University, Princeton, NJ 08540, USA}\\
\textit{e-mail address}: \texttt{gabai@math.princeton.edu}\\

\medskip
\textsc{Department of Algebra and Geometry, Faculty of Science, Palacky University, Olomouc, 771 46 Czech Republic}\\
\textsc{Department of Mathematics, Princeton University, Princeton, NJ 08540, USA}\\
\textit{e-mail address}: \texttt{m.d.trnkova@gmail.com}

\end{document}